\documentclass{amsart}
\usepackage{amsmath,amssymb,bbm,a4,mathrsfs,url,xspace}
\usepackage{array}
\usepackage[all]{xy}
\usepackage{accents}
\usepackage{hyperref}

\newtheorem{Pro}{Proposition}[subsection]
\newtheorem{Le}[Pro]{Lemma}
\newtheorem{Th}[Pro]{Theorem}

\theoremstyle{definition}
\newtheorem{De}[Pro]{Definition}
\theoremstyle{remark}
\newtheorem{Exm}[Pro]{Example}
\newtheorem{Rem}[Pro]{Remark}

\let\x\times
\let\ox\otimes \let\dt\otimes
\let\bt\boxtimes
\def\rbt{\mathop{\accentset{r}{\boxtimes}}\nolimits}
\def\lbt{\mathop{\underaccent{\ell}{\boxtimes}}\nolimits}
\let\ot\otimes

\def\lot{\mathop{\underaccent{\ell}{\otimes}}\nolimits}

\def\blz{\, .\, }

\let\com\diamond
\let\du\vee
\let\bu\bullet

\let\al\alpha
\let\bb\beta
\let\gg\gamma
\let\Gg\Gamma

\let\ee\epsilon

\def\La{\Lambda}
\def\la{\lambda}

\def\A{{\mathbb A}}
\def\B{{\mathbb B}}
\def\D{{\mathbb D}}
\def\G{{\mathbb G}}
\def\K{{\mathbb K}}

\def\sA{{\mathscr A}}
\def\sB{{\mathscr B}}
\def\sC{{\mathscr C}}
\def\sE{{\mathscr E}}
\def\sP{{\mathscr P}}
\def\sV{{\mathscr V}}

\def\f{{\bf f}}
\def\g{{\bf g}}
\def\h{{\bf h}}
\def\i{{\bf i}}
\def\t{{\bf t}}

\def\bhom{{\bf Hom}}
\def\tor{\mathop{\mathrm{Tor}}\nolimits}
\def\cok{\mathop{\mathrm{Coker}}\nolimits}

\def \ker{\mathop{\mathrm{Ker}}\nolimits}
\def \hom{\mathop{\mathrm{Hom}}\nolimits}
\def \End{\mathop{\mathrm{End}}\nolimits}
\def \ext{\mathop{\mathrm{Ext}}\nolimits}
\def\res{\mathop{\mathrm{Res}}\nolimits}
\def\ind{\mathop{\mathrm{Ind}}\nolimits}
\def\H{{\mathrm H}}
\def\E{{\mathrm E}}
\def\es{{\mathfrak S}}
\def\1{^{-1}}

\def\Id{Id}
\def\I{Id}

\let\then\Rightarrow

\begin{document}
\title[Strict polynomial functors and coherent functors]
{Strict Polynomial Functors \\ and coherent functors}
\author[V. Franjou]{Vincent Franjou$^*$}
\address{
Laboratoire Jean Leray \\ Facult\'e des Sciences et Techniques\\
BP 92208\\ F-44322 Nantes Cedex 3}
\email{First.Lastname@univ-nantes.fr}

\thanks{$^{*}$ The first author is partially supported by the LMJL -
Laboratoire de Math\'ematiques Jean Leray, CNRS: Universit\'e de
Nantes, \'Ecole Centrale de Nantes, and acknowledges the hospitality and support of CRM Barcelona where the final corrections on this paper were implemented.}
\author[T. Pirashvili]{Teimuraz Pirashvili}
\address{
Department of Mathematics\\
University of Leicester\\
University Road\\
Leicester\\
LE1 7RH, UK} \email{tp59-at-le.ac.uk}
%
\begin{abstract}
We build an explicit link between coherent functors in the sense of
Auslander \cite{aus} and strict polynomial functors in the sense of
Friedlander and Suslin \cite{FS}. Applications to functor cohomology
are discussed.
\end{abstract}
{\tt\Large Manuscripta DOI: 10.1007/s00229-008-0184-9}
\maketitle
%
\section{Introduction}

Since the foundational work of Schur \cite{schur}, the
representation theory of general linear groups has been closely
related to the representation theory of symmetric groups. Especially
fruitful has been the study, for all integers $n$, of tensor
products $T^n(V):=V^{\otimes n}$ of a vector space $V$, endowed with
the commuting actions of the general linear group $GL(V)$
(diagonally on each factor by linear substitution) and of the
symmetric group $\es_n$ (by permutation of the factors). For many
purposes, the mysterious group ring of the general linear group can
thus be replaced by the more manageable Schur algebra of
$\es_n$-equivariant linear maps of $V^{\otimes n}$.
\par
The use of functors in representation theory, maybe first promoted
by Auslander, is practical and efficient for formalizing the
relations between symmetric groups and general linear groups. The
classical work of Green \cite{green} on representations of the Schur
algebras pushes these ideas quite far and in great generality. In
the example of interest to us, Green associates to every additive
functor $\f$, defined on representations of the symmetric group, the
representation over the Schur algebra given by $\f(V^{\otimes n})$.
The main problem of constructing reverse correspondences is solved
naturally by Green. It is one of the purpose of this paper to shed
new light on these correspondences.
\par\medskip
A few years later, Friedlander and Suslin \cite{FS} introduced
strict polynomial functors, which are equivalent to representations
of the Schur algebra when the dimension of $V$ is at least $n$. This
new formalization is aimed at cohomology computations in positive
characteristic, and it has numerous applications, including a proof
of finite generation of the cohomology of finite group schemes in
the same paper \cite{FS}. The pleasing properties of the category
$\sP$ of strict polynomial functors lead to impressive cohomological
computations \cite{FFSS}, following the fundamental computation of
$\ext_\sP(\I^{(r)},\I^{(r)})$ given in \cite{FS} (the decoration
$(r)$ indicates Frobenius twists, that is extension of scalars through the
$r$-th power of the Frobenius isomorphism). More recently, Cha\l upnik
proves \cite{chal} elegant formulae computing functor cohomology, for various fundamental functors $F$ of the form
$$\f T^n: V\mapsto F(V)=\f(V^{\otimes n}).$$
He succeeds in comparing the groups $\ext_{\sP}(G^{(r)},F^{(r)})=\ext_{\sP}(\f T^{n(r)},\g T^{n(r)})$ with
$\f\, \g\, (\ext_{\sP}(T^{n(r)},T^{n(r)}))$ for many important families
of functors $F$ and $G$. [When the functors $\f$ and $\g$ are given
by idempotents in the group ring of the symmetric group, the two
terms are easily seen to be isomorphic. However, this is rarely
the case in positive characteristic smaller than $n$.] To this end,
Cha\l upnik considers for each strict polynomial functor $F$ certain
adequate choices of functors $\f$, defined on representations of the
symmetric group $\es_n$, such that $\f(V^{\otimes n})=F(V)$, thus
informally rediscovering the correspondences set up by Green. These
methods motivated and inspired the present work.
\par\medskip
Green's correspondences are best expressed in terms of adjoint
functors and recollements of categories. The category $\sP$ thus
appears as a quotient category of the category of all additive
functors defined on representations of the symmetric groups.
Unfortunately, the latter is very large and stays quite mysterious.
Representable functors, such as the functor $\H^0(\es_n,-)$ taking
invariants, are examples of functors obtained through the reverse
correspondences of Green, but they are many more. 
However,  all strict polynomial functors arise from \emph{coherent functors}, which are functors that are presented by representable functors (see
Definition \ref{coherent:De}). The resulting category, if still quite
rich, is much better behaved. For instance, the global dimension of
the category of coherent functors is two. This comes in sharp
contrast with the rich functor cohomology obtained through
homological algebra in the category of strict polynomial functors.
\par
We revisit in this setting some of the properties of functors which
make the category $\sP$ much more tractable than coarse
representations: tensor product, composition (or plethysm),
linearization etc. and we try and find corresponding constructions
for coherent functors. We also apply our insight to functor
cohomology, and obtain Cha\l upnik's constructions in a natural way.

\par\medskip
Section \ref{coherent:section} develops the general properties of
coherent functors. Although we do not claim much originality, it
contains a few results which we could not find in the literature. We
believe Section \ref{tensor of coherent} to be new. 
In Section \ref{category P:section}, we present strict polynomial functors to fit our purpose. 
Section \ref{P vs A:section} contains our main results. It compares coherent functors
and strict polynomial functors. Since the comparison is best stated
in terms of recollements of abelian categories, we recall in an
appendix \ref{FP:appendix} what is needed from this theory. 
Section \ref{tensor:section} lifts the tensor product of polynomial functors
to the level of coherent functors. For later use, the final
section \ref{composition:section} does the same for the composition
of polynomial functors. Section \ref{chal:section} applies this new setting to
functor cohomology and obtains natural versions of Cha\l upnik's
results. Section \ref{chal:section} enjoys its own introduction,
which hopefully makes clear the implications of our results for \cite{chal}, and further developpements.
\par\medskip
{\sc Notations.} We fix a field $\K$ of positive characteristic $p$.
All vector spaces are considered over $\K$, and $\hom$ and $\ox$ are
taken over $\K$, unless otherwise decorated. Let $\sV$ be the
category of finite-dimensional vector spaces. For a finite group
$G$, we let $_G\sV$ denote the category of finite dimensional
$G$-modules.
%

\section{Coherent functors}\label{coherent:section}
Most results in this section are known to experts. A good reference
for the first subsections is a recent survey by Harsthorne
\cite{har}.

\subsection{Finite-dimensional representations as functors}
In this section, we fix a finite group $G$. Let $_G\sV$ be the
category of all finite dimensional $G$-modules and let $\sA(G)$ be
the category of all covariant $\K$-linear functors from $_G\sV$ to
the category $\sV$ of finite dimensional vector spaces. Any $M$ in
$_G\sV$ yields a functor $\t_M$ defined by:
$$\t_M=M\ox_G (-)\ .$$
For instance, the functor $\t_{\K[G]}$ is the forgetful functor to
$\sV$.  Dually, let $\h_M$ be the functor represented by $M$:
$$\h_M=\hom_G(M,-).$$
Note that the functor $\t_{\K[G]}$ is (isomorphic to) the forgetful
functor as well. We shall use that $\h_M(\K[G])$ is isomorphic to
the $\K$-dual vector space $M^{\du}$. This isomorphism is precisely
defined as follows. Let $\tau$ be the element of $ \K[G]^\du$ given
by
$$\tau(g)=0, \ \  {\rm if   } \ \ g\not =1 \ \ \ {\rm and }\ \  \tau(1)=1.$$
By the Yoneda lemma, the function $\tau$ yields a natural morphism
$$\begin{aligned}
\tau_X:\hom_G (X,\K[G])&\to X^\du\\
 f&\mapsto\tau\circ f.
 \end{aligned}$$
The homomorphism $\tau_X$ is an isomorphism when $X=\K[G]$. Since
$\K[G]$ is a self-injective algebra, $\hom_G (-,\K[G])$ is an exact
functor, and $\tau$ is a natural transformation between exact
functors. It results that $\tau_X$ is an isomorphism for all $X$ in
$_G\sV$. We shall use this fact without further reference.
\par
The category $\sA(G)$ is an abelian category. We state below its
elementary properties.
\begin{Pro} For all object $\f$ of the category $\sA(G)$, let $\D (\f)$ be the
object of $\sA(G)$ defined by
$$(\D \f)(M)=(\f(M^\du))^\du .$$
The resulting functor $\D$ is a duality in $\sA(G)$.
\end{Pro}
\begin{Pro} For any $M$ in $_G\sV$, the Yoneda lemma yields a natural
isomorphism
$$\hom_{\sA(G)}(\h_M,\f)\cong \f(M).$$
Thus the functor $\h_M$ is a projective object in the category
$\sA(G)$. Moreover, for all $M,N$ in $_G\sV$, there is a natural
isomorphism:
$$\hom_{\sA(G)}(\h_M,\h_N)\cong \hom_G(N,M).$$
\end{Pro}
\begin{Pro}\label{t:Le} For all $M$ in $_G\sV$, there is a natural isomorphism:
$\D(\h_M)\cong \t_M$. Hence the functor $\t_M$ is an injective
object in the category $\sA(G)$. Moreover, for all $\f$ in $\sA(G)$,
there is a natural isomorphism
$$\hom_{\sA(G)}(\f,\t_M)\cong \D\f(M).$$
In particular, there is a natural isomorphism:
$$\hom_{\sC(G)}(\t_N,\t_M)=\hom_G(N,M).$$
\end{Pro}
\begin{proof}
For $X$ in $\sV$ we define
$$
\theta_X:\hom_G(M,X)\to (X^\du \ox_G M)^\du $$
by:
$\theta_X(\al)=\{\xi\ox m\mapsto\xi(\al(m))\}$. 

It defines a natural transformation $\theta$: $\h_M\to
\D(\t_M)$.
 Since
$\theta_{\K[G]}$ is an isomorphism and both $\h_M$ and $\D(\t _M)$
are left exact functors, it follows that $\theta$ is an isomorphism.
The rest follows because $\D$ is a duality.
\end{proof}
\subsection{Coherent functors}
An object $\f$ in $\sA(G)$ is \emph{finitely generated } if it is a
quotient of $\h_M$ for some $M$ in $\ _G\sV$. Among finitely
generated functors, coherent functors are defined by further
requiring finiteness of the relations.
\begin{De}\label{coherent:De} An object $\f$
in $\sA(G)$ is \emph{coherent}\cite{har}, 
or \emph{finitely presented} \cite[p. 251]{green}\cite[p. 204]{aus},
if it fits in an exact sequence
\begin{equation}\label{carmodgena}
\h_N\to \h_M\to \f\to 0\end{equation} for some $M,N$ in $_G\sV$.
\end{De}
 By the Yoneda lemma, the morphism $\h_N\to \h_M$ is of the form
$\h_\al$ for a uniquely defined $\al:M\to N$. It follows that, for every coherent functor $\f$, there is an exact sequence:
\begin{equation}\label{fisresolventa}
\xymatrix{
0 \ar[r]& \h_{\cok{\al}} \ar[r]& \h_N \ar[r]^{\h_{\al}}& \h_M \ar[r]& \f \ar[r]& 0\ .
}
\end{equation}
 We let $\sC(G)$ be the category of all coherent functors and natural
transformations between them.
It is a classical fact due to Auslander \cite{aus} that the category
$\sC(G)$ is an abelian category with enough projective and injective
objects. Moreover the inclusion $\sC(G)\subset \sA(G)$ is an exact
functor.  
\par
We now give examples of coherent functors.
\begin{Pro}
\begin{enumerate}
    \item For any $M$ in $\ _G\sV$, the functor $\t_M=\D(\h_M)$ is an injective object in
    $\sC(G)$; 
    \item If $\f$ is a
coherent functor, then $\D\f$ is also a coherent functor;
    \item For any integer $i\geq 0$, the homology and
    cohomology functors
$\H ^i(G,-)$, $\H _i(G,-)$ are coherent on $_G\sV$;
    \item For any integer $i\geq 0$, the Tate homology and
cohomology functors $\hat{\H }^i(G,-)$ and $\hat{\H }_i(G,-)$ are
coherent  on $_G\sV$.
\end{enumerate}
\end{Pro}
\begin{proof} Examples (i) to (iii) are also in \cite[\S 2]{har}.
 \begin{enumerate}
    \item If $M$ is free and finite dimensional, there is a natural isomorphism: $\h_{M^\du}\cong\t_M$. 
    For a general $M$, choose a presentation $K\to N\to M\to 0$
in the category $_F\sV$, with free and finite dimensional $K$ and
$N$. Then $\t_M$ is the cokernel of $\t_K\to \t_N$ and hence it is
coherent. It is injective in $\sC(G)$ because it is so even in a
$\sA(G)$ by Proposition \ref{t:Le}.
    \item Assume $\f$ is a coherent functor. By definition, it is a
cokernel of a morphism $\h_N\to \h_M$. Then $\D\f$ is the kernel of
the dual morphism $\t_M\to \t_N$, hence it is also a coherent
functor.
    \item Since
    $$\H ^0(G,-)=\h_\K \hbox{ and } \H _0(G,-)=\t_\K,$$
they are coherent. For a general $i$, choose a projective resolution
$P_*$ of $\K$ with finite dimensional $P_i$, $i\geq0$. Then $\H
^i(G,-)$ is the $i$-th homology of the cochain complex $\h_{P_*}$ of
coherent functors, therefore it is also a coherent functor.
    \item The functors $\hat{\H }^0(G,-)$ and
$\hat{\H }^{-1}(G,-)$ are respectively the cokernel and the kernel
of the norm homomorphism $\H _0(G,-)\to \H ^0(G,-)$, so they are
also coherent. We then proceed as for (iii).
\end{enumerate}
\end{proof}

\subsection{Homological dimension}
For any $M$ in $\ _G\sV$, the functor $\h_M\in \sC(G)$ is
projective, because it is projective even in $\sA(G)$.
The exact sequence (\ref{fisresolventa}) shows that the category $\sC(G)$  has global dimension at most 2. In fact the global dimension is never 1, it is either 2 or 0. The last possibility happens if, and only if,  $\ _G\sV$ is semi-simple \cite{aus} i.e. when the order of $G$ is invertible in $\K$.

The following result characterizes objects of projective or injective dimension smaller than two.
\begin{Pro}\label{80148:Pro} Let $\f$ be a coherent functor.
\begin{enumerate}
    \item The following are equivalent:
    \begin{enumerate}
        \item $\f$ is projective;
        \item The functor $\f$ is left exact;
        \item $\f$ is of the form $\h_M$ for some $M$.
    \end{enumerate}
    \item The following are equivalent:
    \begin{enumerate}
        \item $\f$ is injective;
        \item The functor $\f$ is right exact;
        \item $\f$ is of the form $\t_M$ for some $M$.
    \end{enumerate}
    \item $pd(\f)\leq 1$ if, and only if, $\f$ respect monomorphisms.
    \item $id(\f)\leq 1$ if, and only if, $\f$ respect epimorphisms.
\end{enumerate}
\end{Pro}
\begin{proof}
We prove only (ii) and (iii). The rest follows by duality. Statement
(i) is also proved in \cite[Proposition 3.12 \& 4.9]{har}.
\par
If $\f$ is an injective  object, then $\D(\f)$ is projective, hence
it is a direct summand of an object of the form $\h_M$. It follows
that $\f$ is a direct summand of $\t_M$ for some $M$. The
corresponding projector of $\t_M$ has the form $\t_\al$, where
$\al$ is a projector of $M$. Thus $\f\cong \t_{{\sf Im}(\al)}$. In
particular any injective object is a right exact functor.

\par
Conversely, assume that $\f$ is right exact. The $G$-module
$M=\f(\K[G])$ is finitely generated and therefore one can consider
the functor $\t_M$. There is a well-defined transformation
$\al:\t_M\to \f$ given by: $\al_X(m\ox x)= \f(\hat{x})(m)$, where,
for $x$ in $X$, $\hat{x}:\K[G]\to X$ is the $G$-homomorphism defined
by $\hat{x}(1)=x$. By construction this map is an isomorphism when
$X=\K[G]$. Since $\f$ is right exact $\al_X$ is an isomorphism for
all $X$ in $_G\sV$. Hence $\f$ is injective in $\sC(G)$.
\par
Suppose that $\f$ respects monomorphisms. Consider an exact sequence
of functors
$$0\to \f'\to \h_M\to \f\to 0.$$
We want to show that $\f'$ is projective. By (i), we need to prove
that $\f'$ is left exact. For any short exact sequence in $_G\sV$
$$0\to A\to B\to C\to 0,$$
there is a commutative diagram with exact columns:
$$\xymatrix{&0\ar[d]&0\ar[d]&0\ar[d]&&\\
& \f'(A)\ar[d]\ar[r]&\f'(B)\ar[r]\ar[d]&\f'
(C)\ar[d]\\
0\ar[r]& \h_M(A)\ar[d]\ar[r]&\h_M(B)\ar[r]\ar[d]&\h_M(C)\ar[d]\\
& \f(A)\ar[d]\ar[r]^\al&\f(B)\ar[r]\ar[d]&\f(C)\ar[d]\\
&0&0&0&&}.$$ In this diagram the middle row is also exact and a
diagram chase shows that $\al$ is a monomorphism. It follows that
$\f'$ is left exact and hence projective. This shows that
$pd(\f)\leq 1$.

\par
Conversely, suppose that the projective dimension of $\f$ is $\leq
1$. There is a short exact sequence of functors
$$0\to \h_N\to \h_M\to \f\to 0.$$
If $$0\to A\to B\to C\to 0$$ is a short exact sequence in $_G\sV$,
there is a commutative diagram with exact columns:
$$\xymatrix{&0\ar[d]&0\ar[d]&0\ar[d]&&\\
0\ar[r]& \h_N(A)\ar[d]\ar[r]&\h_N(B)\ar[r]\ar[d]&\h_N(C)\ar[d]\\
0\ar[r]& \h_M(A)\ar[d]\ar[r]&\h_M(B)\ar[r]\ar[d]&\h_M(C)\ar[d]\\
& \f(A)\ar[d]\ar[r]^\al&\f(B)\ar[r]\ar[d]&\f(C)\ar[d]\\
&0&0&0&&}
$$

The first two rows in this diagram are also exact, and $\al$ is a
monomorphism. This shows that $\f$ respects monomorphisms.
\end{proof}

\subsection{Coherent functors and recollement}
 Following Auslander \cite{aus,green}, we relate the category
 of coherent functors $\sC(G)$ with the category $\, _G\sV$. The
 best way to formulate the result is to use the language of
 recollement of categories (see  \cite{kuhn}, \cite{FP},  or Appendix A).

One considers the functors
$$t^*:\,\sC(G)\to \, _G\sV, \  \ \ t_* :  \, _G\sV\to \,\sC(G)
\ \ \  {\rm and } \ \ t_!:\, _G\sV\to \,\sC(G)$$
given respectively by
$$t^*(\f)=\f(\K[G]) , \ \ \ \ t_*(M)=\h_{M^\du}
\ \ {\rm and } \ \ t_!(M):=\t_M.$$
\begin{Pro}\label{ausmicebeba:Le}
\begin{enumerate}
    \item The functor $t_!$ is left adjoint to $t^*$ and, for all $M$ in
$_G\sV$, the $G$-module $t^*t_!(M)$ is naturally isomorphic to $M$.
    \item The functor $\,t_*$ is right  adjoint to
$t^*$ and, for all $M$ in $_G\sV$, the $G$-module $t^*t_*(M)$ is
naturally isomorphic to $M$.
\end{enumerate}
\end{Pro}
\begin{proof}
\begin{enumerate}
    \item For any coherent functor $\f$, the functor $\D\f$ is also
coherent, so one can assume as in (\ref{fisresolventa})
that $\D\f=\cok(\h_{\al}:\h_N\to \h_M)$,
for some linear map $\al$: $M\to N$. 
Applying $\D$, we get an exact sequence:
\begin{equation}\label{fiscoresolution}
\xymatrix{
0 \ar[r]& \f \ar[r]& \t_M  \ar[r]^{\t_{\al}}&  \t_N   .
}
\end{equation}
It shows that:
 $t^*(\f)=\f(\K[G])=
\ker(\al)$. Moreover, there are natural isomorphisms:
$$\begin{aligned}
\hom_{\sC(G)}(\t_M,\f)
&=\hom_{\sC(G)}(\t_M,\ker(\t_\al))\\
&=\ker(\hom_{\sC(G)}(\t_M,\t_\al))\\
&=\ker(\hom_G(M,\al))\\
&=\hom_G(M,\ker\al)\\
&=\hom_G(M,t^*(\f)).
\end{aligned}
$$
From this follows the first statement of (i). The second one follows
from the natural isomorphism: $t^*t_!(M)=\t_M(\K[G])\cong M$.
    \item We use the duality: $\hom_G(M,\K[G])\cong M^\du$.
    Take $\f$ as in the exact sequence (\ref{fisresolventa}):
    $\f=\cok\h_{\al}$ for some linear map $\al$. Because $\K[G]$ is a
self-injective algebra, it follows:
$$t^*(\f)=\f(\K[G])=\cok(\hom_G(\al,\K[G]))\cong
\hom_G(\ker(\al),\K[G])\cong(\ker\, \al)^\du.$$ Moreover, there are
natural isomorphisms:
$$\begin{aligned}
\hom_{\sC}(\f,\h_{X^\du})&= \ker(\h_{X^\du}(\al))\\
&= \ker(\hom_G(X^\du,\al))\\
&= \hom_G(X^\du,\ker\al)\\
&\cong \hom_G((\ker \,\al)^\du, X)\\
&\cong\hom_G(t^*(\f),X)
\end{aligned}
$$
Finally, for $M$ in $ _G\sV$:
$t^*t_*(M)=t^*(\h_{M^\du})=\hom_G(M^\du,\K[G])\cong M$.
\end{enumerate}
\end{proof}

We let ${\sC}^0(G)$ be the full subcategory of ${\sC}(G)$ whose
objects $\f$ are such that: $t^*(\f)=0$. Thus the category
${\sC}^0(G)$ consists exactly of coherent functors which vanish on
projective objects. Since $t^*$ is exact, the subcategory
${\sC}^0(G)$ is abelian. Indeed, it is a Serre subcategory of
${\sC}(G)$. We let
$$r_*:{\sC}^0(G)\to {\sC}(G)$$
be the inclusion. It is an exact functor. It is a consequence of
Proposition \ref{recollements for cheap:Pro} that the functors $r_*,
t^*,t_*,t_!$ are part of a recollement situation
$$
\xymatrix {
 {\sC}^0(G) \ar[rr]|{\ r_*\,}
&&{\sC(G)}\ar@/^3ex/[ll]^{r^!}\ar@/_3ex/[ll]_{r^*}\ar[rr]|{\ t^*\,}
&&{_G\sV}\ar@/^3ex/[ll]^{t_*}\ar@/_3ex/[ll]_{t_!} }.
$$
where $r^*$ and $r^!$, left and right adjoint to $r_*$, are defined
by the following exact sequences:
$$0\to r_*r^!(\f)\to \f\to t_*t^*(\f), \ \ \ t_!t^*(M)\to M\to r_*r^*(\f)\to 0.$$
The following Proposition gives another description of $r^*$ and
$r^!$.
\begin{Pro} \label{801410:Pro}
\begin{enumerate}
    \item For a functor $\f:\, _G\sV\to \sV$, let $\tau:L_0\f\to \f$ be
the natural transformation
 from the 0-th left derived functor. There is an
isomorphism  $$r^*(\f)\cong \cok( \tau).$$
    \item For $M$ in $_G\sV$, let $\Sigma M$ be a finitely generated  $G$-module
which fits in a short exact sequence
\begin{equation}\label{sigmasi}
0\to M\to P\to \Sigma M\to 0\end{equation} where $P$ is a projective
$G$-module. There is an isomorphism
$$r^!(\t _M)\cong \tor^G_1(-, \Sigma M).$$
    \item For a functor $\f:\, _G\sV\to \sV$, write $\f=\ker(\t_\al)$
    for some linear map $\al:M\to N$, as in the exact sequence
(\ref{fiscoresolution}). There is an isomorphism
$$r^!(\f)\cong \ker(\tor^G_1(-, \Sigma \al)).
$$
\end{enumerate}
\end{Pro}
\begin{proof}
\begin{enumerate}
    \item Let us consider a natural transformation
$$\xi:\f\to \g$$
where $\f$ is a coherent functor and $\g$ is in ${\sC}^0(G)$. We
have to prove that $\xi$ factors trough $\cok(\tau)$. In other words
we have to show that
 the composite $\xi\circ \tau:L_0\f\to \g$ is zero.
To this end, for an object $M$ in $\, _G\sV$ choose an exact
sequence $0\to N\to P\to M\to 0$ with projective $P$. The following
commutative diagram with exact top row implies the result:
$$\xymatrix{L_0\f(P)\ar[r]\ar[d]_{\cong}&
L_0\f(M)\ar[d]\ar[r]& 0\\
 \f(P)\ar[r]\ar[d]&
\f(M)\ar[d]&\\
0=\g (P)\ar[r]&\g(M)&}.$$
    \item The long exact sequence for
$\tor$-groups on $0\to M\to P\to \Sigma M\to 0$ yields an exact
sequence
$$0\to \tor^G_1(-, \Sigma M)\to \t_M\to \t_P\to \t_{\Sigma M}\to
0.$$
Let $\psi:\t_M\to \h_{M^\du}$ be the natural transformation
defined by:
$$\begin{aligned}
\psi_X: X&\ox _GM&\to&& \hom_G(M^\du,X)\\
x&\ox m&\mapsto&& (\xi\mapsto\xi(m)x).
 \end{aligned}$$
It is an isomorphism when $M$ is a projective object in $_G\sV$.
There is a commutative diagram with exact rows
$$\xymatrix{&0\ar[r]& \h_{M^\du}\ar[r]&
\h_{P^\du}\ar[r]&\h_{\Sigma M^\du}&\\
0\ar[r]&\tor^G_1(-,\Sigma_M)\ar[r] &
\t_M\ar[r]\ar[u]&\t_P\ar[r]\ar[u]_\cong &\t_{\Sigma
M}\ar[r]\ar[u]&0.}$$ It follows that there is an exact sequence
$$0\to \tor^G_1(-, \Sigma M)\to \t _M\to \h_{M^\du}.$$
The result follows from the comparison with the exact sequence
$$0\to
r_*r^!(\t_M)\to \t_M\to t_*t^*(\t_M))=t_*(M)\cong \h_{M^\du}.$$
    \item Apply $r^!$ to the
exact sequence (\ref{fiscoresolution}). Because the functor $r^!$ is
left exact:
$$r^!(\f)=r^!(\ker(\t_\al)=\ker(r^!(\t_\al)),$$
and the result follows from (ii).
\end{enumerate}
\end{proof}

\subsection{Induction and restriction for coherent functors}
Let $H$ be a subgroup of the finite group $G$. The induction functor
$\ind ^G_H:{\, _H\sV}\to {\, _G\sV}$ is left and right adjoint
to the restriction functor $\res ^G_H:{\, _G\sV}\to {\, _H\sV}$.
Both functors are exact. By pre-composition, we obtain exact
functors
$$\uparrow^G_H:\sA(H)\to \sA(G)\ \ \text{and}\ \
\downarrow^G_H:\sA(G)\to \sA(H)$$ which are given by
$$(\uparrow^G_H \g)(M)=\g(\res _H^G(M))\ \ \text{and}\ \
(\downarrow^G_H \f)(N)=\f(\ind _H^G(N))$$ for $M\in \, _G\sV$,
$N\in \, _H\sV$, $\f\in \sA(G)$ and $\g\in \sA(H)$. From the adjunction of induction and restriction, it is formal
that $\downarrow^G_H$ is left and right adjoint to the functor
$\uparrow^G_H$.

\begin{Pro} If $\f$ is coherent in $\sC(G)$, then $\downarrow^G_H(\f)$
is coherent in $\sC(H)$. If $\g$ is coherent in $\sC(H)$, then
$\uparrow^G_H(\g)$ is coherent in $\sC(G)$. Thus, induction and
restriction define adjoint pairs of functors
$$\uparrow^G_H:\sC(H)\to \sC(G)\ \ \text{and}\ \
\downarrow^G_H:\sC(G)\to \sC(H).$$
\end{Pro}

\begin{proof} Take $M\in \, _G\sV$ and $N\in \, _H\sV$. Since $\ind $ is left adjoint to $\res $ we have
$$\uparrow^G_H(\h_N)\cong \h_{\ind ^G_H N}$$
Similarly, since $\ind $ is also right adjoint to $\res $:
$$ \downarrow^G_H(\h_M)=\h_{\res ^G_HM}.$$
The rest follows from the exactness of $\uparrow$ and $\downarrow$
and the exact sequence (\ref{carmodgena}).
\end{proof}

\subsection{Universal property of the category of coherent functors}
The following elementary proposition stresses the relevance of the
category of coherent functors.

\begin{Pro}\label{cuni} Let $\sE$ be an abelian category and assume additive functors
$U:\, _G\sV\to \sE$ and  $V:\, (_G\sV)^{op}\to \sE$ are given. Then
\begin{enumerate}
    \item  There exists a left
exact functor $U^c:\sC(G)\to \sE$, unique up to an isomorphism,
 such that $U^c(\t_\al)=U(\al)$ for all $\al$ in $\sC(G)$. 
 Moreover $U^c$ is exact when $U$ is right exact.
 \item  There exists a
 right exact functor $V^c:\sC(G)\to \sE$, unique up to an isomorphism
 such that $V^c(\h_\al)=V(\al)$ for all $\al$ in $\sC(G)$. 
 Moreover $V^c$ is exact when $V$ is left exact.

\item  Suppose that $U$ is right exact, that $V$ is left exact and
suppose that natural isomorphisms $\al_P:U(P)\to V(P^\du)$ are
given for all  projective objects $P$ of the category $\, _G\sV$. Then there exists an isomorphism of functors
$\al^c:U^c\to V^c$ such that for any projective object $P\in \,
_G\sV$ the morphism $\al^c(\t_P)$ is the same as the following
composite:
$$U^c(\t_P)=U(P)\to V(P^\du)=V^c(\h_{P^\du})\cong V^c(\t_P),$$
where the last isomorphism is induced by the isomorphism
$\t_P\cong \h_{P^\du}$.
\end{enumerate}
\end{Pro}
\begin{proof} 
Deriving the additive functor provides the desired extension. The details are as follows.
\par
i) Let $\f$ be a coherent functor. Write $\f=\ker\t_{\al}$ as in (\ref{fiscoresolution}) to get an injective resolution of $\f$:
\[\xymatrix{
0 \ar[r]& \f \ar[r]& \t_M  \ar[r]^{\t_{\al}}&  \t_N  \ar[r]&  \t_{\cok(\al)} \ar[r]& 0\ .
}\]
Thus $U^c(\f)$ is necessarily the kernel of $U(\al):U(M)\to U(N)$.
Not only this proves the uniqueness of the functor $U^c$, but it
also gives the construction of $U^c$ as the $0$-th right derived
functor of the restriction of $U^c$ on injective objects. If $U$ is
right exact then
$$0\to U^c(\f)\to U(M)\to U(N)\to U(\cok(\al))\to 0$$
is an exact sequence. Thus the right 1-st (and therefore any higher)
derived functor vanishes and hence $U^c$ is an exact functor.
Similarly for ii).

iii) For $M\in \, _G\sV$ choose a free presentation, an exact sequence
$$ K\to L\to M\to 0$$
with free $K$ and $L$ in  $\, _G\sV$. 
This gives an exact sequence of coherent functors
$$\t_K\to \t_L\to \t_M\to 0.$$
Since both functors $U^c$ and $V^c$ are exact, there exists a unique
morphism $\al^c(\t_M)$ which fits in a commutative diagram with
$\al^c(\t_K)$ and $\al^c(\t_L)$.
  Since any coherent functor $\f$ has a resolution
$0\to\f\to \t_M\to \t_N$ and $U^c$, $V^c$ are exact functors,
there exists a unique morphism $\al^c(\f)$ which fits in a
commutative diagram with $\al^c(\t_M)$ and $\al^c(\t_N)$. The
claim follows.
\end{proof}

\subsection{External tensor products of coherent functors}\label{tensor of coherent}
Let $G$ and $H$ be finite groups. For $M$ in $\, _G\sV$ and $N$ in
$\, _H\sV$, let $M\boxtimes N$ denote the vector space $M\ox N$ with
$G\x H$-action. It yields an exact bifunctor
$$\boxtimes: \, _G\sV\x _H\sV\to \, _{G\x H}\sV.$$
We wish to extend this bifunctor to coherent functors. 

We first discuss a concrete elementary way of extending the product, inspired by \cite[p. 781]{chal}.
Note that for $X$ in $_{G\x H}\sV$ and $\g$ in $\sA(H)$, 
the vector space $\g (\res ^{G\x H}_{1\x H}(X))$ has a natural $G$-structure, because $G\x 1$ commutes with
$1\x H$ in $G\x H$. Thus, for $\f$ in $\sA(G)$, the evaluation
$\f(\g(\res ^{G\x H}_{1\x H} X))$ makes sense. 
Denote it by $(\f \odot\g)(X)$. We have thus defined an additive bifunctor
$$\odot: \sA(G)\x \sA(H)\to \sA(G\x H).$$
\begin{Pro}\label{odot:Pro}
Let $G$ and $H$ be finite groups.
\begin{enumerate}
 \item For all coherent $\f$ in $\sC(G)$ and $\g$ in $\sC(H)$, the functor
 $\f\odot \g$ is coherent in $\sC(G\x H)$. This defines a bifunctor
\[
\odot: \sC(G)\x \sC(H)\to \sC(G\x H).
\]
\item For all $M$ in $_{G}\sV$ and $N$ in $_{H}\sV$, there are natural isomorphisms:
\[
\h_M\odot \h _N\cong  \h_{M\bt N}\,,\ \ \t_M\odot\t _N\cong  \t_{M\bt N}.
\]
\item The bifunctor $\odot$ is exact with respect to the first
argument.
\item The bifunctor $\odot$ commutes with duality.
\end{enumerate}
\end{Pro}
\begin{proof} The last two points are easy. We prove the first two.
\begin{enumerate}
\item[(ii)] For $X$ in $_{G\x H}\sV$, one has:
\[
\begin{aligned}
\h_{M\boxtimes N}(X)&=
\hom_{G\x H}(M\boxtimes N,X)\\
&= \hom_{G}(M, \hom_{H}(N,\res ^{G\x H}_{1\x H}X)\\
&=\h_M\odot \h_N(X).
\end{aligned}
\]
and the other isomorphism follows by duality (or directly as easily).
\item[(i)] By (iii), it is enough to prove that $\h_{M}\odot\g$ is coherent for all $\g$ in $\sA(H)$ and  $M$ in $_{G}\sV$. Taking a free presentation of $M$ reduces further to the case when $M$ is free. Because $\h_{M}\odot\g$ is then exact in $\g$, the case $\g=\h_{N}$ is enough. This is (ii).
\end{enumerate}
\end{proof}
\begin{Exm}\label{odot:Exm}
Note that: $\f\odot\h_{\K}(X)=\f(X^H)$ but: $\h_{\K}\odot\f(X)=\f(X)^H$, so that $\h_{\K}\odot\f$ is not even right exact in $\f$, even for $H=\es_{2}$.
\end{Exm}
As a direct application of the method of Proposition \ref{cuni}, we obtain symmetric exact alternatives to this extension of the external tensor.
\begin{Pro}\label{lt:Pro}
There is a right exact (in each argument) functor:
\[
\lbt: \sC(G)\x\sC(H)\to\sC(G\x H),
\]
unique up to isomorphism, such that: 
\[
\h_M\lbt\h_N=\h_{M\bt N}
\]
for all $M$ in $\, _G\sV$ and $N$ in $\, _H\sV$.
This bifunctor is symmetric and balanced.
\end{Pro}
\begin{proof}
For $M$ in $\sV_{H\x G}$, let $M^\flat$ be the $G\x
H$-module with same underlying vector space and action given by:
$(g,h)._{\flat}m=(h,g). m$. 
Symmetry means that the following diagram commutes up to
natural isomorphism:
   $$\xymatrix{\sC(G)\x \sC(H)\ar[r]^{\boxtimes}\ar[d]^{tw}&\sC(G\x
   H)\ar[d]^{Tw}\\
   \sC(H)\x \sC(G)\ar[r]_{\boxtimes}&\sC(H\x G)}$$
where $tw$ is twisting of factors, while $Tw$ is given by
$Tw(\f)(M)=\f(M^{\flat})$, for $\f$ in $\sC(G\x H)$.
To prove this property, it suffices to note that both functors $Tw\circ \lbt$ and $\lbt \circ tw$ are right exact and take isomorphic
values on projective generators; for projective generators, it reduces to the symmetry of the external tensor product.

To prove that the bifunctor is balanced, it is enough, by symmetry, to prove that $-\lbt\h_{N}$ is exact for each $N$ in $_{H}\sV$. This functor is the extension to $\sC(G)$ of the left exact functor $M\mapsto \h_{M\bt N}$. It is therefore exact by Proposition \ref{cuni}(ii).
\end{proof}We now compare these constructions.
\begin{Pro} \label{lttoodot:Pro}
There is a natural transformation: 
\[
\f\lbt\g\to\f\odot\g
\]
for all $\f$ in $\sC(G)$ and $\g$ in $\sC(H)$, which is an isomorphism when $\g$ is projective.
\end{Pro}
\begin{proof}
Since both sides are right exact on $\f$, the isomorphism for a projective generator $\g=\h_{N}$ follows from the case of projective generators $\f=\h_M$. 
\end{proof}
\begin{Exm}\label{odot vs. lt:Exm}
Note that the two products do not always coincide. For instance, using Example \ref{odot:Exm} for $G=H$, 
$\h_{\K}\lbt\f\cong\f\lbt\h_{\K}\cong\f\odot\h_{\K}$ is not always equal to $\h_{\K}\odot\f$.
\end{Exm}
\begin{Pro}\label{rt:Pro}
There is a left exact (in each argument) functor:
\[
\rbt: \sC(G)\x\sC(H)\to\sC(G\x H),
\]
unique up to isomorphism, such that: 
\[
\t_M\rbt\t_N=\t_{M\bt N}
\]
for all $M$ in $\, _G\sV$ and $N$ in $\, _H\sV$.
This bifunctor is symmetric and balanced. 
\end{Pro}
\begin{Pro}\label{odottort:Pro}
There is a natural transformation: 
\[
\f\odot\g\to\f\rbt\g
\]
for all $\f$ in $\sC(G)$ and $\g$ in $\sC(H)$, which is an isomorphism when $\g$ is injective.
\end{Pro}
These can be deduced from the natural isomorphism:
\[
\D(\f\lbt\g)\cong\D\f\rbt\D\g.
\]
\begin{Rem}\label{()=rt:Rem}
Indeed, the $\rbt$-product coincide with Cha\l upnik's mysterious $(-,-)$ construction \cite[p. 785]{chal}.
For, choose an embedding of $\f$ in some $\f'$ and of $\g$ in some $\g'$.
Because the $\rbt$-product is left exact, $\f\rbt\g$
embeds in $\f'\rbt\g'$. More precisely, it is contained both in
the image of $\f'\rbt\g$ and in the image of $\f\rbt\g'$
in $\f'\rbt\g'$. A simple diagram chase shows that
$\f\rbt\g$ indeed \emph{equals the intersection of these
images} in $\f'\rbt\g'$.
Choosing $\f'$ and $\g'$ to be injectives gives a description in term of the more expicit $\odot$-product and recovers Cha\l upnik's definition of $(\f,\g)$ as well as the definition of $\rbt$ as a $0$-th right derived functor of $\odot$.
\end{Rem}
\subsection{Extension of the domain of coherent functors}\label{extactcoh}
In this section, we show that a coherent functor $\f$ in $\sC(G)$
can be evaluated on any object in an abelian category which is
equipped with a $G$-action.

Let $\sE$ be a $\K$-linear abelian category. There exists a unique (up to a unique isomorphism) exact bifunctor
$$\ox: \sV\x \sE\to \sE$$
such that $\K\ox(-):\sE\to \sE$ is isomorphic to the identity functor
$\Id_\sE:\sE\to \sE$. It is obtained as follows. 
On the full subcategory of $\sV$ with objects
$\K^n$, $n\geq 0$, it is defined by: $\K^n\ox E:=E^n$, naturally in
$E$. We are left with defining the functor $(-)\ox E:\sV\to \sE$ on
morphisms, for a fixed object $E$. For a matrix $(k_{ij})$,
$k_{ij}\in \K$, the corresponding morphism $E^n\to E^m$ is given by
the same matrix, but where $k_{ij}$ is considered as an elements of
$\End_{\sE}(E)$ via the $K$-linear structure on $\sE$.

Putting $$\bhom(V,A):=V^\du\ox A$$ defines an  exact bifunctor
$$\bhom: (\sV)^{op}\x \sE\to \sE$$ such that $\bhom(\K,-):\sE\to
\sE$ is isomorphic to the identity functor $\Id_\sE:\sE\to \sE$. It
is clear that one has a natural isomorphism
$$\bhom(V,W\ox A)\cong W\ox \bhom(V,A),$$
for $V,W$ in $\sV$ and $A$ in $\sE$. Thus $\sE$ is tensored and
co-tensored over $\sV$.
\begin{Rem}\label{structural hom iso}
More generally, for any $\K$-linear functor $T:\sE\to\sE'$ between
$\K$-linear categories, there are natural isomorphisms: $V\ox
T(A)\cong T(V\ox A)$ and $\bhom(V,TA)\cong T(\bhom(V,A))$ for $V$ in
$\sV$ and $A$ in $\sE$. Indeed, the second isomorphism is a
consequence of the first one, while the first one is obvious for
$V=\K$ and follows from the additivity of $T$ for all $V$'s.
\end{Rem}
We are now ready to define invariants and co-invariants of a finite
group $G$ in an abelian category $\sE$. A \emph{left $G$-object in
$\sE$} is a pair $(A,\al)$, where $A$ is in $\sE$ and $\al:G\to {\sf
Aut}_{\sE}(A)$ is a group homomorphism. Equivalently it is an object
$A$ equipped with a map $\la:\K[G]\ox A\to A$ satisfying obvious
properties. We let ${\sf Rep}(G,\sE)$ be the abelian category of
left $G$-objects in $\sE$. If $A$ is a $G$-object, its structural
map $\K[G]\ox A\to A$ has an adjoint $A\to \bhom(\K[G],A)$. One
defines:
$$\H^0(G,A):=\ker(A\to \bhom(\K[G],A)).$$
If $M$ is a $G$-module, one defines $\bhom_G(M,A)$ to be
$\H^0(G,\bhom(M,A))$.
\par
Similarly, the reader will define $M\ox _G A$ as a coequalizer of two
canonical maps $M\ox \K[G]\ox A\to M\ox A$. We also put:
$\H_0(G,A):=\K\ox_G A$.

\begin{Pro} There exists a bifunctor
$$<-,->\,:\sC(G)\x {\sf Rep}(G,\sE)\to \sE$$
unique up to an isomorphism, such that the functor $<-,A>:\sC(G)\to
\sE$ is exact for any $G$-object $A$ and
$$<\t_M,A>\, =\, M\ox_GA$$
for any $M\in \, _G\sV$. Moreover, one has a natural isomorphism
$$<\h_M, A>\, \cong \, \bhom_G(M,A).$$
\end{Pro}

\begin{proof} Fix a $G$-object $A$. Put $R(M)=M\ox_GA$ in Proposition \ref{cuni}
to get the existence and uniqueness of the pairing. To show the last
assertion, observe that both sides of the isomorphism are left exact
additive functors on $M$ and therefore it suffices to prove the
isomorphism for $M=\K[G]$. Since $\h_{\K[G]}\cong \t_{\K[G]}$ is
the forgetful functor, both sides in this case are isomorphic to
$A$.
\end{proof}

 It is clear that for $\sE=\,_G\sV$ this pairing is just evaluation:
$<\f,M>=\f(M).$ The pairing allows to evaluate a coherent functor on
graded $G$-objects, for instance.

\section{Application to strict polynomial functors}
\subsection{The category of strict polynomial functors}\label{category P:section}
Strict polynomial functors were introduced by Friedlander and Suslin
\cite[\S 2]{FS}. We now explain what they are in a given degree, in
a way suitable for our purpose (see also \cite[\S 4]{pira}).

Let us fix a positive integer $n$. We start with the $n$-th divided
power of a vector space $V$, defined by:
$$\Gg^n(V):=\H^0(\es_n, V^{\ox n})=(V^{\ox n})^{\es_n},$$
where the symmetric group on $n$-letters $\es_n$ acts
 on $V^{\ox n}$ by permuting the factors.
For any $x$ in some vector space $X$, we let $\gg(x)$ be the element
$x^{\ox n}$ in $\Gg^n(X)$. This defines a natural set map $\gg_X:X\to
\Gg^n(X)$. Reordering the factors
$$A^{\ox n}\ox B^{\ox n}\to (A\ox B)^{\ox n}$$
induces a $\K$-linear natural transformation
$$\Gg ^n(A)\ox \Gg^n(B)\to \Gg^n(A\ox B)$$
sending $\gg_A(a)\ox \gg_B(b)$ to $\gg_{A\ox B}(a\ox b)$. Together with
the composition law in $\sV$, these maps define a composition map:
$$
\Gg^n(\hom(V,W))\ox \Gg^n(\hom(U,V))\to
\Gg^n(\hom(V,W)\ox\hom(U,V))\to\Gg^n(\hom(U,W)).
$$
This defines a category $\Gg^n\sV $, with the same objects as $\sV$,
and with morphisms
$$\hom_{\Gg^n\sV}(V,W):=\Gg^n(\hom(V,W)).$$
The following Lemma describes the category $\Gg^n\sV$ as a full
subcategory of $_{\es_n}\sV$.
\begin{Le}\label{i:Le} For a positive integer $n$, the functor
$$\begin{aligned}
\i:\Gg^n\sV&\to \, _{\es_n}\sV\\
    V&\mapsto V^{\ox n}
    \end{aligned}$$
is a full embedding.
\end{Le}
\begin{proof}
This follows from the natural isomorphism:
$$\hom_{\es_n}(V^{\ox n}, W^{\ox n})=
(\hom(V^{\ox n},W^{\ox n}))^{\es _n}\cong (\hom(V,W)^{\ox n})^{\es _n}=
\hom_{\Gg^n\sV}(V,W).$$
\end{proof}
Reformulating \cite[\S 2]{FS}, a \emph{homogeneous strict polynomial
functor of degree $n$} defined on $\sV$ is a $\K$-linear functor
$\Gg^n\sV\to \sV$. We let $\sP_n$ be the category of homogeneous
strict polynomial functors of degree $n$. It is known \cite[\S
3]{FS} that the category $\sP_n$ is equivalent to the category of
finite dimensional modules over the Schur algebra $S(n,n)$.
\par
The collection of maps $\gg_X:X\to \Gg^n(X)$ yields a (nonlinear)
functor $\gg:\sV\to \Gg^n\sV$. Pre-composition with $\gg$ associates
to any strict polynomial functor defined on $\sV$ an usual functor
on $\sV$; it is called the underlying functor of the strict
polynomial functor. It is usual to denote by the same letter a
strict polynomial functor and its underlying functor. For example,
the composite
$$\xymatrix{\Gg^n\sV\ar[r]^\i & \  _{\es_n}\sV \ar[rr]^{H^0(\es _n,-)}&& \sV}$$
is denoted by $\Gg^n$, since its underlying functor is the $n$-th
divided power functor, and $S^n$ denotes the composite
$$\xymatrix{\Gg^n\sV\ar[r]^\i & \  _{\es_n}\sV \ar[rr]^{\H _0(\es_n,-)}&& \sV},$$
because its underlying functor is  the $n$-th symmetric  power.
Similarly the composite
$$\xymatrix{\Gg^n\sV\ar[r]^\i & \  _{\es_n}\sV \ar[rr]^{\text{forget}}& &\sV}$$
 is denoted by $T^n$,
because the underlying functor is the $n$-th tensor power. We now
recall from \cite[\S 3]{FS} the basic properties of the category
$\sP_n$.
\par
There is a well-defined tensor product of strict
polynomial functors which corresponds to the usual tensor product of
underlying functors, and it yields a bifunctor
$$-\ox -:\sP_n\x \sP_m\to \sP_{n+m}.$$
For example: $T^n=T^1\ox \cdots \ox T^1$ ($n$-factors).
\par
There is also a duality in $\sP_n$. For an object $F$ in $\sP_n$, we
let $\D F$ be the homogeneous strict polynomial functor given by
$$(\D F)(V)=(F(V^\du))^\du$$
where $W^\du$ denotes the dual vector space of $W$. Since the values
of any  homogeneous strict polynomial functor are finite
dimensional, $\D$ is an involution and defines an equivalence of
categories $\D:\sP_n^{op}\to \sP_n$. The functor $\D F$ is called
the \emph{dual} of $F$.
\par
 The category $\sP_n$ has enough projective
and injective objects. A set of generators is indexed by partitions
of  $n$, that is decreasing sequences of positive integers adding up
to $n$. For a partition $\la=(n_1\geq n_2\geq \cdots \geq n_k)$, put
$$\Gg^\la:=\Gg^{n_1}\ox \cdots \ox \Gg^{n_k}.$$
The functors $\Gg^\la$, when $\la$ runs through all partitions of
the integer $n$, form a set of projective generators. Indeed,
$\hom_{\sP_n}(\Gg^\la,F)$ is the evaluation on the base field of the
cross-effect of the functor $F$ of homogeneous multidegree $\la$
(i.e. the component of weight $\la$ under the action of $\G_{m}^k$ on the cross-effect of $F$).
Dually, the functors $$S^\la:=S^{n_1}\ox \cdots \ox S^{n_k} $$ form a
set of injective cogenerators. In particular, the functor
 $T^n$ is projective and injective in $\sP_n$. Moreover, the action
of $\es_n$ by permuting factors yields an exact functor
$$\begin{aligned}
c^*:\sP_n\to&\, _{\es_n}\, \sV \\
F\mapsto &\hom_{\sP_n}(T^n,F).
\end{aligned}$$
The representation $c^*(F)$ is often called the \emph{linearization}
of the functor $F$; we use the letter $c$ for cross-effect.
The functor $c^*$ has both a left and a right adjoint functor given
respectively by
$$(c_!(M))(V)=(M\ox V^{\ox n})_{\es_n}\ , \ \
c_*(M)=(M\ox V^{\ox n})^{\es_n}.$$ Let $\sP^0_n$ be the full
subcategory of $\sP_n$ whose objects are the strict polynomial
functors $F$ such that $c^*(F)=0$. This condition means that the
underlying functor has degree less than $n$ in the additive sense of
Eilenberg and MacLane. Let $d_*:\sP_n^0\to \sP_n$ be the inclusion
and let $d^*$ and $d^!$ be the left and right adjoint of $d_*$. By
Proposition \ref{recollements for cheap:Pro}, this defines a
recollement situation:
$$
\xymatrix {
 {\sP_n}^0 \ar[rr]|{\ d_*\,}
&&{\sP_n}\ar@/^3ex/[ll]^{d^!}\ar@/_3ex/[ll]_{d^*}\ar[rr]|{\ c^*\,}
&&{\, _{\es_n}\sV}\ar@/^3ex/[ll]^{c_*}\ar@/_3ex/[ll]_{c_!} }.
$$
More on this recollement can be found in \cite{kuhnstrat}.

\subsection{The relation between  strict polynomial and coherent functors}\label{P vs A:section}
For simplicity, we write $\sA_n$ and  $\sC_n$ instead of
$\sA(\es_n)$ and $\sC(\es_n)$. For an object $\f$ in $\sA_n$, the
composite
$$\xymatrix{\Gg^n\sV\ar[r]^\i & \  _{\es_n}\sV \ar[r] ^\f& \sV}$$
defines a strict polynomial functor, which is denoted by $j^*(\f)$.
This construction defines a functor
$$\begin{aligned}
&j^*:\sA_n\to \sP_n\\
&j^*(\f)(V)=\f(V^{\ox d}). \end{aligned}$$ The same functor was
constructed, in terms of Schur algebras, by Green \cite[\S 5, pp
275--276]{green}. The functor $j^*$ was also considered by Cha\l
upnik \cite[\S 3]{chal} in his work on functor cohomology.
\begin{Le} The functor $j^*$ respects duality:
$\D\circ j^*\cong j^*\circ \D$.
\end{Le}
\begin{proof}
$\D j^*(\f)(V)=(j^*(\f)(V^\du))^\du=(\f(V^{\du \ox n}))^\du \cong (\D
\f)(V^{\ox n})=j^*\D(\f)(V).$
\end{proof}
\begin{Pro} The functor $j^*:\sA_n\to \sP_n$ has a right adjoint
functor defined by:
$$j_*(F)(M)=\hom_{\sP_n}(j^*(\h_M),F)=
\hom_{\sP_n}(\hom_{\es_n}(M,(-)^{\ox n}),F)$$ where $M$ is
representation of $\es_n$. It has also a left adjoint functor
defined by:
$$j_!(F)=\D( j_*(\D F)).$$
In other words:
$$(j_!F)(M)^\du=\hom_{\sP_n}(F, (-)^{\ox n}\ox_{\es_n}M^\du).$$
\end{Pro}
\begin{proof}Since hom's in the category $\sP_n$ are finite dimensional
vector spaces, we see that $j_*(F)$ belongs to $\sA_n$. The fact
that it is right adjoint of $j^*$ follows from the Yoneda lemma. The
dual formula is formal:
$$\hom_{\sA_n}(\D (j_*(\D F)),\f)
\cong \hom_{\sA_n}(\D\f, j_*(\D F))\cong \hom_{\sP_n}(j^*(\D \f),\D F)$$
$$\cong \hom_{\sP_n}(\D j^*(\f),\D F)\cong \hom_{\sP_n}(F, j^*(\f)).$$
To check the last formula, observe that:
$$j_*(\D F)(M^\du)=\hom_{\sP_n}(j^*(\h_{M^\du}),\D F)\cong \hom_{\sP_n}(F,j^*(\D
\h_{M^\du}))$$ and
$$
j^*(\D \h_{M^\du}) \cong j^*\t_{M^\du}=(-)^{\ox n}\ox_{\es_n}M^\du.
$$
\end{proof}
\begin{Rem}  In particular $j_*$ and $j_!$ are a functorial
choice of, respectively, an injective and projective symmetrization
of \cite[Section 3]{chal}.
\end{Rem}
\begin{Rem}
The existence of adjoints of a precomposition functor is quite a
general phenomenon, see Example \ref{recollements for functors:Exm}.
\end{Rem}
We now study these adjoint functors. For a partition $\la$ of a
positive integer $n$, we let $\es_\la$ be the corresponding Young
subgroup of $\es_n$.
\begin{Le}\label{sym:Le} Let $\la$ be a partition  of $n$ there are natural
isomorphisms:
$$j_*(S^\la)\cong \H _0(\es_\la,-)=\t_{\K[\es_n/\es_\la]},$$
$$j_!(\Gg^\la)\cong \H ^0(\es_\la,-)=\h_{\K[\es_n/\es_\la]}.$$
In particular
$$\hom_{\sP_n}(S^\la,S^\mu)=\hom_{\es_n}(\K[\es_n/\es_\la],
\K[\es_n/\es_\mu])=\hom_{\sP_n}(\Gg^\mu,\Gg^\la).$$
\end{Le}
\begin{proof}
$$\begin{aligned}
j_*(S^\la)(M)
&=\hom_{\sP_n}(\hom_{\es_n}(M,(-)^{\ox n}),S^\la)\\
&\cong \hom_{\sP_n}(\Gg^\la, \D(\hom_{\es_n}(M,(-)^{\ox n})))\\
&\cong \hom_{\sP_n}(\Gg^\la, (-)^{\ox n}\ox_{\es_n}M).
\end{aligned}$$
The first isomorphism follows from \cite[Corollary 2.12]{FS}. The
second follows by duality. The last statement follows from the fact
that $j_!$ is full and faithful.
\end{proof}
\begin{Pro}\cite[2.14]{green} 
The values of the functors $j_*$ and $j_!$ are coherent functors.
\end{Pro}
\begin{proof}
By duality, it is enough to consider the functor $j_*$. By Lemma
\ref{sym:Le}, the statement is true for injective cogenerators of
$\sP_n$. Since $j_*$ is left exact the result follows by a
resolution argument.
\end{proof}
\begin{Pro}
The unit $\Id_{\sP_n}\to j^*j_!$ and the counit $j^*j_*\to
\Id_{\sP_n}$ are isomorphisms.
\end{Pro}
\begin{proof} We prove only the second isomorphism, the first one follows by
duality. It is clear that $j^*(\H _0(\es_\la,-))\cong S^\la$. Thus
Lemma \ref{sym:Le} shows that the statement is true for injective
cogenerators of $\sP_n$. Since $j^*$ is exact and $j_*$ is left
exact, the result follows by taking resolutions.
\end{proof}
\begin{Rem}
Since the functor  $j_*$ is a full embedding, this gives a new proof
of \cite[Lemma 3.4]{chal}.
\end{Rem}
\begin{Pro}\label{recollement coherent polynomial:Pro}
Let ${\sC}_n^{\sf Y}$ be the full subcategory of ${\sC}_n$ whose
objects are the coherent functors $\f$ such that, for all partitions
$\lambda$ of $n$:
$$\f(\K[\es_n/\es_\lambda])=0.$$
The functors $j^*$ and its adjoints $j_*,j_!$ are part of a
recollement of abelian categories:
$$
\xymatrix {
 {\sC}_n^{\sf Y}\ar[rr]|{\ i_*\,}
&&\sC _n \ar@/^3ex/[ll]^{i^!}\ar@/_3ex/[ll]_{i^*}\ar[rr]|{\ j^*\,}
&&\sP_n\ar@/^3ex/[ll]^{j_*}\ar@/_3ex/[ll]_{j_!} }
$$
\end{Pro}
\begin{proof}
According to Proposition \ref{recollements for cheap:Pro}, the
functor $j^*$ and its adjoints give rise to a recollement situation.
To determine the kernel category, it is enough to notice that every
Young-permutation representation $\K[\es_n/\es_\lambda]$ occurs as a
direct factor in the tensor product $V^{\ox n}$ as soon as the
dimension of $V$ is $n$.
\end{proof}
\begin{Exm}  ${\sC}_n^{\sf Y}=0$ if $p>n$. Moreover ${\sC}_n^{\sf Y}=0$ for $p=2$ and $n=2$ or $3$.\end{Exm}
\begin{Pro}\label{09314:Pro}
The counit $j_!j^*(\f)\to \f$ is an isomorphism when $\f=\t_M$.
Dually, the unit $\f\to j_*j^*(\f)$ is an isomorphism when
$\f=\h_M$.
\end{Pro}
\begin{proof} We prove only the first assertion. Since both
$j_!j^*(\t_M)$ and $\t_M$ are right exact functors of $M$, it is
enough to consider the case when $M=\K[\es_n]$. In this case $\t_M$
is the forgetful functor. Therefore: $j^*(\t_M)=\ox
^{n}=\Gg^{11\cdots 1}$ and
$$j_!j^*(\t_{\K[\es_n]})=j_!(\Gg^{11\cdots
1})=\H ^0(\es_{11\cdots 1},-)=\t_{\K[\es_n]}.
$$
\end{proof}
\begin{Pro}\label{norm iso:Pro} The norm transformation  (see Appendix \ref{FP:appendix})
 for the previous recollement
situation is an isomorphism on projective and injective objects.
\end{Pro}
\begin{proof}
By Lemma \ref{sym:Le} we have $j_!(\Gg^\lambda)=\h_M$, for
$M=\K[\es_n/\es_\lambda]$. By Proposition \ref{09314:Pro} we have
also $j_*(\Gg^\lambda)=j_*j^*(\h_M)=\h_M$, thus the norm is an
isomorphism on projective objects. By duality it is also an
isomorphism on injective objects.
\end{proof}
The following examples gather some other known values of the adjoints $j_!$, $j_*$.
\begin{Exm} The relations of Proposition \ref{09314:Pro}:
 $$j_*(j^*\t_M)\cong\t_M\cong j_!(j^*\t_M)$$
apply in particular when $M$ is the signature, or when $M$ is
induced from the signature of a Young subgroup $\es_\mu$. This shows
that, when $p$ is odd, the norm is an isomorphism on a tensor
product of exterior powers $\La^\mu$.
\end{Exm}
\begin{Exm} Let $\I^{(1)}$ be the Frobenius twist in $\sP_p$,
that is extension of scalars through the Frobenius \cite[p. 212]{FS}. It is related to
the norm $N$ by an exact sequence:
$$\xymatrix{0\ar[r]&\I^{(1)}\ar[r]&S^p\ar[r]^N&\Gg^p\ar[r]&\I^{(1)}\ar[r]&0}$$
It follows that there are exact sequences:
$$0\to j_*(\I^{(1)})\to \H _0(\es_p,-)\to \H ^0(\es_p,-)$$
$$ \H _0(\es_p,-)\to \H ^0(\es_p,-) \to j_!(\I^{(1)})\to 0.$$
Thus $j_*(\I^{(1)})=\hat{\H }^{-1}(\es_p,-)$ and
$j_!(\I^{(1)})=\hat{\H }^0(\es_p,-)$ (where as usual $\hat{\H }$
denotes Tate cohomology).
\end{Exm}
\begin{Exm} Assume $p=2$. Let $S$ be a set with $n$ elements. For each $0\leq
k\leq n$ we let $B_k$ be the vector space spanned on the set of all
subsets of $S$ with exactly $k$-elements. Define $d:B_k\to B_{k+1}$
by
$$d(X)=\sum_{X\subset Y\in B_{k+1}}Y$$
Then $d^2=0$ and $B_*$ is a cochain complex of $\es_n$-modules. One
checks that it is acyclic. For an integer $m\geq 1$ and
$n=2^{m+1}$, the explicit injective resolution of $S^{2^m(1)}$ of
\cite[\S 2]{FLS} and \cite[\S 8]{FS} allows to compute:
\begin{equation}\label{chal1401}
R^*j_*(S^{2^m(1)})\cong \H ^*(\t_{B_*}).\end{equation} In
particular, $j_*(S^{2^m(1)})$ is the kernel of the obvious map $\H
_0(\es_n,-)\to \H _0(\es_{n-1,1},-)$. Another consequence of
(\ref{chal1401}) is the fact that $R^kj_*(S^{2^m(1)})=0$ when $k\geq
m$.
\end{Exm}
We now show the compatibility of the different recollement
situations.
\begin{Pro} There are commutative diagrams of categories and functors:
$$\xymatrix{ &\\&\\
 \sC_n&\sP_n\ar[l]_{j_*}\\
 _{\es_n}\sV\ar@{=}[r]\ar[u]^{t_*}&_{\es_n}\sV\ar[u]^{c_*}
}
\xymatrix{\sC_n^{\sf Y}\ar[r]\ar@{=}[d]&\sC_n^{0}\ar[r]
\ar[d]_{r_*}&\sP^0_n\ar[d]^{d_*}\\
\sC_n^{\sf Y}\ar[r]^{i_*}& \sC_n\ar[r]^{j^*}\ar[d]_{t^*}&
\sP_n\ar[d]_{c^*}\\
& _{\es_n}\sV\ar@{=}[r]&_{\es_n}\sV } \ \ \ \ \ \ \xymatrix{&\\
 \sC_n&\sP_n\ar[l]_{j_!}\\
 _{\es_n}\sV\ar@{=}[r]\ar[u]^{t_!}&_{\es_n}\sV\ar[u]^{c_!}.
}$$
\end{Pro}
\begin{proof}
To show that $c^*\circ j^*=t^*$, note that the three functors
involved are exact. It is therefore enough to check that they
coincide on $\h_M$ for each $M$ in $_{\es_n}\sV$. This means that we
need to show that:
$$\hom_{\sP_n}(T^n,\hom_{\es_n}(M,(-)^{\ox n}))\cong M^\du$$
Since $T^n$ is projective, the left hand side of the expected
isomorphism is left exact as a functor of $M$. So it suffices to
consider the case when $M$ is injective, and it reduces to the case
when $M=\K[\es_n]$. In this case it is a well-known isomorphism:
$$\hom_{\sP_n}(T^n,T^n)\cong \K[\es_n].$$
\par
To show that $j_!c!=t!$, note that both sides are right exact. It is
therefore enough to check that they coincide on $\K[\es_n]$. In this
case, $j_!c_!(\K[\es_n])=j_!(T^n)$ has already been seen (see the
proof of Proposition \ref{09314:Pro}) to be the forgetful functor
$t_{\K[\es_n]}$.\par The rest is quite similar.
\end{proof}


%

\subsection{Tensor products of coherent functors}\label{tensor:section}
The aim of this section is to lift the bifunctor $\ox:\sP_n\x\sP_m\to
\sP_{n+m}$ given by $(F\ox G)(V)=F(V)\ox G(V)$ at the level of
coherent functors. Not surprisingly, it involves the induction functors
\[
\ind ^{\es_{m+n}}_{\es_m\x \es_n}:\, _{\es_m\x \es_n}\sV\to\,
_{\es_{n+m}}\sV .
\]
For $M$ in $_{\es_m}\sV$and $N$ in $_{\es_n}\sV$, we let $M\boxtimes N$
denote the vector space $M\ox N$ with $\es_m\x \es_n$-action, and
simply denote by $M\ox N$ the module $\ind ^{\es_{m+n}}_{\es_m\x
\es_n}(M\boxtimes N)$. The operation $\otimes$ yields a symmetric
monoidal structure on the category $\bigoplus_{n\geq 0}{\ _{\es_n}\sV}$.

Using Section \ref{tensor of coherent}, one defines biadditive functors  
$\sC_m\x \sC_n\to \sC_{m+n}$ by:
\[
\f \lot \g= \ind^{\es_{m+n}}_{\es_n\x \es_m}(\f\lbt \g) \ \ \ \
\f \blz \g= \ind^{\es_{m+n}}_{\es_n\x \es_m}(\f\odot\g).
\]
\begin{Th} \label{tensor in coherent:Th}
\begin{enumerate}
    \item For all $M$ in $_{\es_m}\sV$ and $N$ in $_{\es_n}\sV$, there
 are natural isomorphisms
 $$\h_M\lot\h _N\cong\h_{M\otimes N}= \h_M\blz\h_N.$$
\item For all $\f$ in $\sC_m$ and all  $\g$ in $\sC_{n}$,
there is a natural transformation:
\[
\f \lot \g \to \f \blz \g
\]
which is an isomorphism when $\g$ is projective.
\item The bifunctor $\lot$ equips the category $\bigoplus _{d\geq 0} \sC_d$
with a symmetric monoidal structure.
\item The bifunctor $\lot$ is right exact (in each argument) and balanced.
 \item For all $\f$ in $\sC_m$ and $\g$ in $\sC_n$, there are natural
isomorphisms
$$j^*(\f)\ox j^*(\g)\cong j^*(\f\lot\g)\cong j^*(\f\blz\g);$$
    \item For all $F$ in $\sP_m$ and $G$ in $\sP_n$, there is a natural
isomorphism
$$j_!(F)\lot j_!(G)\cong j_!(F\ox G)$$
\end{enumerate}
\end{Th}
\begin{proof} Right exactness of $\lot$ and (i) follow from Proposition \ref{lttoodot:Pro}. As in Proposition \ref{cuni}, these properties characterize $\lot$ up to isomorphism, and (iii) follows.
We next show the first isomorphism in (v). Since
$j^*$ and $\lot$ are right exact functors, it is enough to consider the
case when $\f=\h_M$ and $\g=\h_N$.  We need to prove that
$$\h_M(V^{\ox m})\ox \h_N(V^{\ox n})\cong\h_{M\otimes N}(V^{\ox m+n}).$$
This follows from the isomorphisms:
\[\begin{aligned}
\hom_{\es_m}(M,V^{\ox m}) \ox \hom_{\es_n}(N,V^{\ox n})
&= \hom_{\es_m\x \es_n}(M\bt N,V^{\ox m}\bt V^{\ox n})\\
&=\hom_{\es_m\x \es_n}(M\bt N,\res _{\es_m \x \es_n}^{\es_{m+n}}V^{\ox m+n})\\
&=\hom_{\es_{m+ n}}(M\ot N,V^{\ox m+n}).
\end{aligned}\]
To show the second isomorphism in (v), one observes that 
for all $W\in\sV$ and $N\in\, _{\es_n}\sV$, there is
a natural isomorphism (as in Remark \ref{structural hom iso}):
$${\bf g}(W\ox N)\cong W\ox {\bf g}(N),$$
where $\es_n$ acts on the second factor of $W\ox N$. Now: $V^{\ox
m+n}=V^{\ox m}\ox V^{\ox n}$ as $\es_n$-modules, with trivial action on
the first factor, so:
$$\g( \res ^{\es_{m+n}}_{1_m\x \es_n} V^{\ox m+n})=
V^{\ox m}\ox\g(V^{\ox n}).$$ Similarly, the group $\es_m$ acts
trivially on $\g(V^{\ox n})$, and we obtain: 
$$
j^*(\f\blz \g)(V)=\f(\g(\res ^{\es_{m+n}}_{1_m\x \es_n}V^{\ox m+n}))\cong\f(V^{\ox m})\ox \g(V^{\ox n}).
$$
Finally, we show the isomorphism in (vi). Let us denote by $\la \cup \mu$ the
concatenation of two partitions $\la$ and $\mu$. There is an isomorphism:
\begin{equation}\label{cup:equation}
\ind ^{\es_m}_{\es _\la}\K\ox \ind ^{\es_n}_{\es _\mu}\K =
{\ind}^{\es_{m+m}}_{\es_{\la \cup \mu}}\K.
\end{equation}
Now to (vi). It is enough to assume that $F$ and $G$ are projective
generators: $F=\Gg^\la$ and $G= \Gg^\mu$. By Lemma
\ref{sym:Le}, we have
$$j_!(\Gg^\la)=\H ^0(\es_\la,-)=\h_{\K[\es_m/\es_\la]}=\h_{\ind ^{\es_m}_{\es _\la}\K}\ , \ \
j_!(\Gg^\mu)=\h_{\ind ^{\es_n}_{\es _\mu}\K}.
$$
So the isomorphism (vi) in this case follows from the isomorphism
(\ref{cup:equation}). 
\end{proof}
We leave to the reader the statement dual to Theorem \ref{tensor in coherent:Th}.

\section{Application to functor cohomology}\label{chal:section}
\subsection{Introduction}
The first computation of $\ext$-groups between strict polynomial
functors in \cite{FS} states that the graded vector space
$A_r=\ext^*_{\sP_{p^r}}(\I^{(r)},\I^{(r)})$  is one-dimensional in
even degrees smaller than $2p^r$ and zero else. A lot more basic
computations were carried out in \cite[\S 4 - 5]{FFSS}. Based on
these results, M. Cha\l upnik \cite[\S 4 - 5]{chal} has succeeded in
extending them to include a lot more basic functors. His
calculations eventually rely on the very special form of the
fundamental but elementary computation of
$\ext^*_{\sP_{np^r}}(\I^{(r)\ox n},\I^{(r)\ox n})$ as $B_r= A_r^{\ox d}\ox
\K[\es_d]$, a $\es_n^{op}\x\es_n$-permutation graded module.

The present paper's setting allows to better formulate the tools
Cha\l upnik used and the results he obtained. This includes his
notion of symmetrization \cite[\S 3]{chal}, which we proved in
Section \ref{P vs A:section} to be functorial, or his
$(-,-)$-product \cite[Theorem 4.4, p. 785]{chal}, which coincide
with our $\rbt$-product of Section \ref{tensor of coherent} by
Remark \ref{()=rt:Rem}. We thus explain key results showing
exactness and symmetry in the contravariant/covariant variables, for
example. We apply our insight to prove naturality of \cite[Theorem
5.3]{chal}, a result which does not follow from the methods in
\cite{chal}, and we further extend Cha\l upnik's results. It has to
be noted however, that the results in \cite{chal} do not follow
formally from our considerations, but they rather are the elementary
calculations to build on from.
\subsection{Cohomology of Frobenius twists}
We start (as does Cha\l upnik \cite{chal}) by recalling
\cite[Theorem 4.5]{FFSS}.

Let us observe that in the category $\sP:=\bigoplus_{d\geq 0}\sP_d$,
a tensor product of two projectives is projective. This allows a
K\"unneth morphism:
$$
\ext_\sP^*(A,B)\ox \ext_\sP^*(C,D)\to\ext_\sP^*(A\ox C,B\ox D).
$$
There is thus a natural map:
$$
\ext_\sP^*(A,B)^{\ox d}\to\ext_\sP^*(A^{\ox d},B^{\ox d})
$$
which defines a natural homomorphism
$$
S^d(\ext^*_{\sP}(A,B))\to \ext^*_{\sP}(\Gg^d\circ A,S^d\circ B).
$$
Taking $A=B=\I^{(r)}$, the $r$-th Frobenius twist of the identity
functor, we get the map
\begin{equation}\label{smap}
S^d(\ext^*_{\sP_{p^r}}(\I^{(r)},\I^{(r)}))\to
\ext^*_{\sP_{dp^r}}(\Gg^{d(r)},S^{d(r)}) \end{equation} which, by a
special case of Theorem 4.5 in \cite{FFSS}, is an isomorphism.

From this, Cha\l upnik \cite[Theorem 4.3]{chal} deduces that, when
$G=\Gg^\mu$, there is a natural in $G$ isomorphism:
\begin{equation}\label{odot}
\ext_{\sP}^*(G^{(r)},F^{(r)})\ \cong \ <j_*F \odot  j_*\D G, B_r>.
 \end{equation}
Our reformulation shows naturality in $F$ as well. Note that:
$j_*\D G=\t_{\K[\es_d/\es_\mu]}$, so by Proposition \ref{odottort:Pro}:
\begin{equation}\label{boxtimes:eq}
j_*F \odot  j_*\D G=j_*F \rbt j_*\D G.
 \end{equation}
On the latter formulation, left exactness in both $F$ and $G$, as
well as symmetry, become transparent.
\par
 We apply the formula for $G=\Gg^d$, knowing the $\es_d$-bimodule structure on $B_r$.
We readily obtain the following natural version of \cite[Theorem
5.3]{chal}:
\begin{Pro}
For all $F$ in $\sP_d$, there is a natural isomorphism:
$$
\ext_{\sP_{dp^i}}(\Gg^{d(i)},F^{(i)})\cong
F(\ext_{\sP_{p^i}}(\I^{(i)},\I^{(i)})).
$$
\end{Pro}
More generally:
\begin{Pro} Let $F$ and $G$ be strict polynomial functors of
degree $d$. There exists a natural spectral sequence:
 $$E^{i, *}_2=\ <R^iT(G,F), B_r> \, \then \,
 \ext_{\sP}^{i+*} (G^{(r)},F^{(r)}),$$
 where $R^iT$ denotes the $i$-th right derived functor of the left exact
 functor 
 \[\begin{aligned}
 T:\sP_d^{op}\x \sP_d&\to \sC(\es_d^{op}\x\es_{d})\\
(G,F) &\mapsto j_*\D G\rbt j_*F.
 \end{aligned}\]
\end{Pro}
\begin{proof} First consider the special case when $G$ is
projective. In this case we have to show that
\[
<j_*\D G\rbt j_*F, B_r> \, \cong \, \ext_{\sP}^* (G^{(r)},F^{(r)})
\]
Since both sides are additive in $G$, it is enough to check the
isomorphism when $G=\Gg^\mu$. This case results from (\ref{odot})
and (\ref{boxtimes:eq}).

For the general case, we take a projective resolution of $G_\bu$ of
$G$. After twisting, we obtain a (non-projective) resolution
$G^{(r)}_\bu$ of $G^{(r)}$. The hypercohomology spectral sequence
obtained by applying $\hom_{\sP}(-,F^{(r)})$ to $G^{(r)}_\bu$ has
the form
\[
E^{i,j}_1=\ext^j_{\sP}(G^{(r)}_i,F^{(r)})\, \then \,
\ext_{\sP}^{i+j} (G^{(r)},F^{(r)}).
\]
Since $G_i$ is projective we can apply the previous computation to obtain:
\[
E^{i,*}_1=\ <j_*\D G\rbt j_*F_i , B_r> 
\]
and the $E_2$-term has the expected form. Left exactness of the $\rbt$-product ensures that in the first column appears $R^{0}T=T$.
\end{proof}

\appendix
\section{}
\subsection{Recollement of abelian categories}\label{FP:appendix}
To reveal the relationship between the different abelian categories,
we use the language of recollements (see \cite{kuhn} and \cite{FP}).
A recollement of abelian categories consists of a diagram of abelian
categories and additive functors
$$
\xymatrix {
 {\sA}'\ar[rr]|{\ i_*\,}
&&{\sA} \ar@/^3ex/[ll]^{i^!}\ar@/_3ex/[ll]_{i^*}\ar[rr]|{\ j^*\,}
&&{\sA}''\ar@/^3ex/[ll]^{j_*}\ar@/_3ex/[ll]_{j_!} }
$$
satisfying the following conditions:
\begin{enumerate}
    \item the functor $j_!$ is left adjoint to $j^*$ and the functor $j^*$ is left adjoint of
$j_*$;
    \item the unit  $Id_{{\sA}''}\to j^*j_!$ and the counit $j^*j_*\to
  Id_{{\sA}''}$ are isomorphisms;
    \item the functor $i^*$ is left adjoint of $i_*$ and $i_*$ is left adjoint of
$i^!$;
    \item the unit  $Id_{{\sA}'}\to i^!i_*$ and the counit  $i^*i_*\to
Id_{{\sA}'}$ are isomorphisms;
    \item the functor $i_*: {\sA}'\to {\ker }(j^*)$ is an equivalence of
categories.
\end{enumerate}
\begin{Exm}
The following example is the paradigm of a recollement situation.
Let $X$ be a space, $C$ is a closed subset in $X$ and $U=X\setminus
C$ its open complement. Extension and restriction yield a
recollement of sheaves categories:
$$
\xymatrix {
 {\sf Sh}(C)\ar[rr]|{i_*}
&&{\sf Sh} (X) \ar@/^3ex/[ll]^{i^!}\ar@/_3ex/[ll]_{i^*}\ar[rr]|{\
j^*\,} &&{\sf Sh}(U)\ar@/^3ex/[ll]^{j_*}\ar@/_3ex/[ll]_{j_!} }.
$$
\end{Exm}
\par
The list of properties (i)-(v) can be somewhat shortened.
\begin{Pro}\label{recollements for cheap:Pro}
Let  $j^*:{\sA }\to {\sA }''$ be an exact functor which satisfies
{\rm (i)} and {\rm (ii)}: it admits both a left adjoint $j_!$ and a
right adjoint $j_*$, and the unit $Id_{{\sA }''}\to j^*j_!$ and
counit $j^*j_*\to Id_{{\sA }''}$ are isomorphisms. Let ${\sA }'$ be
the full subcategory of $\sA $ with objects those $A$ such that
$j^*A=0$. Then the full embedding $i_*:{\sA }'\to {\sA }$ has
adjoint functors $(i^*,i^!)$ and the unit $Id_{{\sA }'}\to i^!i_*$
and counit $i^*i_*\to Id_{{\sA }'}$ are isomorphisms. In other words
we have a recollement situation.
\end{Pro}
\begin{proof} Let $A$ in $\sA$ and let $\epsilon_A$: $j_!j^*A\to A$ be the
counit of the adjoint pair $(j_!,j^*)$. Because $Id_{\sA ''}\to
j^*j_!$ is an isomorphism, we have $j^*(\cok \ee_A)=0$. It follows
that ${\cok}(\ee_A)$ lies in the subcategory $\sA '$. So there is a
well-defined functor $i^*:{\sA }\to {\sA '}$ such that
${\cok}(\epsilon_A)=i_*i^*A$. The rest follows, using the short
exact sequence of natural transformations:
$$ j_!j^* \buildrel \epsilon \over\rightarrow Id_{\sA }\to i_*i^*\to 0$$
and the dual study of the unit of adjunction $\eta$ which sums up in
the following exact sequence:
$$0\to i_*i^!\to Id_{\sA } \buildrel \eta \over\rightarrow j_*j^* \ .$$
\end{proof}
\begin{Rem} Actually, if $\sA $ is a category
 of modules over a ring (or, more generally, if $\sA $
 is a Grothendieck category),
 then it is enough to assume that $j^*$ is an exact functor which
 has a left adjoint functor $j_!$ such that the unit of
 adjunction:
 $Id_{{\sA }''}\to j^*j_!$
 is an isomorphism. The existence of $j_*$ follows from
 \cite[Proposition 2.2]{nt}.
\end{Rem}
\begin{Exm}\label{recollements for functors:Exm}
Recollements arise naturally when relating functor categories
through precomposition. Indeed, starting with a functor
$\i:\sA\to\sB$, precomposition is an exact functor:
$$\begin{aligned}
j^*:\ \sV^\sB&\to\sV^\sA\\
F&\mapsto F\circ \i .
\end{aligned}$$
A classic result of D. Kan tells that it always admits adjoint
functors, called the left and the right Kan extension. By \cite[\S
X.3, Corollary 3]{macl}, the unit and the counit of adjunction are
isomorphisms when the functor $\i$ is a full embedding. A
recollement situation then arises by Proposition \ref{recollements
for cheap:Pro}.
\par
In the case when $\i$ is a full embedding of $\K$-linear categories,
the functor $j^*$ and its adjoints restrict to the subcategories of
$\K$-linear functors. Proposition \ref{recollement coherent
polynomial:Pro} describes the resulting recollement when the functor
$\i$ is the full embedding of Lemma \ref{i:Le}.
\end{Exm}
Another useful functor arises from a recollement: the functor
$j_{!*}:\sA''\to\sA$ is the image of the norm $N:j_!\to j_*$, the
natural transformation which
 corresponds to $1_X$, under the isomorphism
$$
\hom_{\sA}(j_!X,j_*X)\cong \hom_{\sA''}(X,j^*j_*X)\cong
\hom_{\sA''}(X,X).$$ The functor $j_{!*}$ preserves simple objects,
and every simple object in $\sA$ is either the image of a simple in
$\sA'$ by the functor $i_*$, or the image of a simple in $\sA''$ by
the functor $j_{!*}$.
\par
We close with an immediate consequence of the Grothendieck spectral
sequence for a composite functor.
\begin{Pro}\label{sse:Pro} Assume in a recollement situation all abelian categories have
enough projective objects. For $X$ in $\sA''$ and $B$ in $\sA$,
there are spectral sequences:
$$\E^{pq}_2=\ext^p_{\sA}(L_qj_!(X),B)\Longrightarrow \ext^{p+q}_{\sA''}(X,j^*B)$$
and
$$\E^{pq}_2=\ext^p_{\sA}(B, R^qj_*(X))\Longrightarrow \ext^{p+q}_{\sA''}(j^*B,X).$$
\end{Pro}

\subsection{Composition and coherent functors}\label{composition:section}
The composite of two strict polynomial functors is a strict polynomial functor
 \cite{FS}. For completness we  lift the resulting bifunctor:
$$\begin{aligned}
\circ:\sP_n\x \sP_m&\to \sP_{nm}\\
(F,G)&\mapsto F\circ G \end{aligned}
$$
at the level of coherent functors.
\par

Composition of functors is
exact with respect to the first variable. Although the functor
$G\mapsto F\circ G$ is not additive for $n>1$, it still has some
exactness properties.
\begin{De} Let $\A$ and $\B$ be abelian categories.
For any short exact sequence
\begin{equation}\label{ses:equation}
\xymatrix{0 \ar[r]&A_1\ar[r]^\al& A\ar[r]^{\bb}&A_2\ar[r]&0}
\end{equation} in $\A$, define $\delta_1, \delta_2: A\oplus A_1\to
A$ by
$$\delta_1(a,a_1)=a+\al(a_1)\ , \ \ \ \delta_2(a,a_1)=a\ .$$
A covariant functor $T:\A\to
\B$ \emph{preserves reflective coequalizers} if for every short exact sequence (\ref{ses:equation}), the following
sequence is exact: 
\begin{equation}\label{T(ses):equation}
\xymatrix{
T(A\oplus A_1)\ar[rrr]^{T(\delta_1)-T(\delta_2)}&&&T(A)\ar[r]^{T(\bb)}&T(A_2)\ar[r]&0
}.
\end{equation}
\end{De}
Observe that when $T$ is  additive, then it preserves reflective
coequalizers if, and only if, $T$ is  right exact.  Let us observe
also that if the exact sequence (\ref{ses:equation}) splits then the
sequence (\ref{T(ses):equation}) is exact for any functor $T$. If
$\A$ has enough projective objects, then any (possibly nonadditive)
functor, from the category of projective objects in $\A$ to the
category $\B$, has a unique (up to unique isomorphism) extension to
a functor $\A\to \B$ which preserves reflective coequalizers.

\begin{Le} For any $F$ in $\sP_n$, the functor
$$\begin{aligned}
\sP_m\to &\sP_{nm}\\
G\mapsto &F\circ G \end{aligned}$$
 preserves reflective
coequalizers and coreflective equalizers.
\end{Le}
\begin{proof} Take any short exact sequence in $\sP_m$:
$$\xymatrix{0\ar[r]& G_1\ar[r]^\al& G\ar[r]^{\bb}&G_2\ar[r]&0}.$$
After evaluating at $V\in \sV$, the corresponding sequence
$$0\to G_1(V)\to G(V)\to G_2(V)\to 0$$
splits. Therefore for any $F$, the sequence
$$F(G(V)\oplus G_1(V))\to F(G(V))\to F(G_2(V))\to 0$$
is exact. This shows that $F\circ (-)$ respects reflective
coequalizers. Similarly for coreflective equalizers.
\end{proof}
For a natural number $m$ and a group $G$, let $\es_m\wr G$ be the
wreath product, which by definition is the semi-direct product
$G^m\ltimes \es_m$. For $M$ in $_{\es_m}\sV$ and $N$ in $_{G}\sV$,
it acts on $M\ox N^{\ox m}$. In particular, for $G=\es_n$, let
$$M\bullet N:= \ind  ^{\es_{mn}}_{\es_m\wr\es_n}(M\ox N^{\ox m})$$
It defines a functor:
$$\bullet:{ _{\es_m}\sV}\x {\, _{\es_n}\sV}\to\, _{\es_{mn}}\sV$$
\begin{Pro} There is a unique (up to isomorphism) bifunctor
$$\com:\sC_m\x \sC_n\to \sC_{mn}$$
with the following properties
\begin{enumerate}
    \item  The functor $\com$ is exact with respect of the first variable and
 preserves reflective
    coequalizers with respect to the second variable;
    \item  $\h_M\com \h_N=\h_{M\bullet N}$.
\end{enumerate}
\end{Pro}

We define another bifunctor
$$\bar{\circ}:\sC_m\x \sC_n\to \sC_{mn}.$$
Because $\dt$ is a symmetric monoidal structure on $\oplus_{d\geq
0}\sC_d$, the functor $\g^{\dt m}$ has a natural action of $\es_m$
for $\g$ in $\sC_n$. We put (compare with \cite{benoit}):
$$\f\bar{\circ} \g:=<\f,\g^{\dt m}>.$$
\begin{Pro} For all $\f$ in $\sC_m$ and $N$ in $_{\es _n}\sV$, there is a natural
isomorphism:
$$\f\com \h_N\cong \f\bar{\circ}\h_N.$$
\end{Pro}
\begin{proof} Since $\f\com \g$  and
 $ \f\bar{\circ}\g$ are exact on $\f$, it is enough to consider the case
$\f=\h_M$. Then:
$$
\begin{aligned}
\h_{M\bullet N}(X)&=\hom_{\es_{mn}}(\ind ^{\es_{mn}}_{\es_n\wr \es_m}
(M\ox N^{\ox m}),X)\\
&=\hom_{\es_n\wr \es_m}(M\ox N^{\ox m},\res ^{\es_{mn}}_{\es_n\wr \es_m}(
X))\\
& =\hom_{\es_m}(M,  \hom_{(\es_n)^m}(N^{\ox m},\res ^{\es_{mn}}_
{(\es_n)^m}(X)))\\
&=\h_M(\h_{N^{\otimes m}})(X) \\
&=(\h_m\bar{\circ}\h_N)(X).
\end{aligned}
$$
\end{proof}
\begin{Pro}
\begin{enumerate}
 \item For all $\f$ in $\sC_m$ and $\g$ in $\sC_n$, there is a natural
isomorphism
$$j^*(\f)\circ j^*(\g)\cong j^*(\f\com \g).$$
    \item For all $F$ in $\sP_m$ and $G$ in $\sP_n$, there is a natural
isomorphism
$$j_!(F)\com j_!(G)\cong j_!(F\circ G).$$

\item
for all $F$ in $\sP_m$ and $G$ in $\sP_n$, there is a natural
isomorphism
$$j_*(F)\com \, j_*(G)\cong j_*(F\circ G),$$

\item $\t_M\com\t_N=\t_{M\bullet N}$;

\end{enumerate}
\end{Pro}
\begin{proof}
\begin{enumerate}
    \item Since $j^*$ is right exact and
 $\com$ respects reflective coequalizers,
it suffices to consider the case when $\f=\h_M$ and $\g=\h_N$. Then
we have:
$$
\begin{aligned} j^*(\f)\circ j^*(\g)&\cong \hom_{\es_m}(M,
(\hom_{es_n}(N,V^{\ox n}))^{\ox m})\\
&=\hom_{\es_m}(M,  \hom_{(\es_n)^m}(N^{\ox m},\res ^{\es_{mn}}_
{(\es_n)^m}(V^{\ox mn})))\\
\end{aligned}
$$
By the previous computation the last group is isomorphic to
$\h_{M\bullet N}(V^{\ox mn})$.
    \item It is enough to assume
that $F$ and $G$ are projective generators: $F=\Gg^{\mu}$,
$G=\Gg^{\nu}$. We set $M=\K[\es_m/\es_\mu]$ and
$N=\K[\es_n/\es_{\nu}]$. Then we have $j_!(F)=\h_M$ and
$j_!(T)=\h_N$. Therefore $j_!(F)\com j_!(G)\cong \h_{M\bullet N}$.
On the other hand
$$\begin{aligned}
F\circ G (V)&=\Gg^mu(\Gg^\nu(V))
\hom_{\es_m}(M,(\hom_{\es_n}(N,V^{\ox n}))^{\ox m})\\
&=\hom_{\es_m}(M,\hom_{\es_n\x \cdots \x \es_n}(N^{\ox m}\res ^{\K[\es_{mn}]}_{\es_n\x\cdots\x\es_n}(V^{\ox mn})))\\
&=\hom_{\es_m\wr \es_n}(M\ox N^{\ox m}, \res ^{\K[\es_{mn}]}_{es_m\wr \es_n}(V^{mn}))\\
&=\hom_{\es_{mn}}(N\bullet M, V^{\ox mn})
\end{aligned}
$$
It follows that: $j_!(F\circ G)=\h_{M\bullet N}$.
\end{enumerate}
\end{proof}



\begin{thebibliography}{33}

\bibitem{aus} M. Auslander. {\it Coherent functors}. Proc. Conf.
Categorical Algebra (La Jolla, Calif., 1965)  Springer (1966), 189--231.

\bibitem{chal} M. Cha\l upnik. {\it  Extensions of strict
polynomial functors.}  Ann. Sci. \'Ecole Norm. Sup. (4) 38 (2005),
no. 5, 773--792.




\bibitem{FP} V. Franjou \& T. Pirashvili.
{\it Comparison of abelian categories recollements}.  Doc. Math.  9  (2004), 41--56.

\bibitem{benoit} B. Fresse. {\it On the homotopy of simplicial algebras
over an operad.}  Trans. Amer. Math. Soc.  352  (2000),  no. 9,
4113--4141.

\bibitem{FLS} {V.  Franjou, J. Lannes \& L. Schwartz.} 
{\it Autour de la cohomologie de MacLane des corps finis}.
 Invent. Math.  115  (1994),  no. 3, 513--538. 

\bibitem{FFSS} {V.  Franjou, E. Friedlander, A. Scorichenko \& A.
Suslin.} {\it General linear and functor cohomology over finite
fields}. Ann. of Math. (2) 150 (1999), 663--728.

\bibitem{FS} E. Friedlander \& A. Suslin.{\it Cohomology of finite
group schemes over a field}.  Invent. Math.  127  (1997),  no. 2, 209--270.

\bibitem{green} J.A. Green. {\it  On three functors of M.
Auslander's.}
Comm. Algebra 15 (1987), 241--277.

\bibitem{har} R. Harsthorne. {\it Coherent functors}. Adv. in Math. 140
(1998),  no. 1, 44--94.

\bibitem{kuhn}  N. J. Kuhn.
{\it Generic representations of the finite general linear groups and
the Steenrod algebra}. II. $K$-Theory 8 (1994), no. 4, 395--428.

\bibitem{kuhnstrat} N. J. Kuhn.
{\it A stratification of generic representation theory and
generalized Schur algebras}. $K$-Theory 26 (2002), no. 1, 15--49.

\bibitem{macl} S. MacLane. Categories for the Working
Mathematician. Graduate Texts in Mathematics 5, Springer, New York-Berlin, 1971. ix+262 pp.


\bibitem{nt} {C. Nastasescu \& B. Torrecillas}. {\it Colocalization on
Grothendieck categories with applications to coalgebras}. J.
Algebra, 185 (1996), 108--124.

\bibitem{pira}{T. Pirashvili}. {\it Introduction to functor cohomology.}
\emph{In} Rational Representations, the Steenrod Algebra, and Functor Homology, 
Panor. Synth\`eses {16}, Soc. Math. France, Paris, 2003.

\bibitem{schur}  I. Schur. Thesis (1903). \emph{In} Gesammelte Abhandlungen, Band I, 5-76. Springer, 1973.

\end{thebibliography}
\end{document}